\documentclass[reqno]{amsart}
\usepackage{color}

\numberwithin{equation}{section}
\usepackage{latexsym}
\usepackage{upgreek}
\usepackage{amsmath}
\usepackage{amssymb}
\usepackage{mathrsfs}
\usepackage{graphicx}
\usepackage{ifthen}
\usepackage[T1]{fontenc}
\usepackage[latin1]{inputenc}

\numberwithin{equation}{section}

\newtheorem{defi}{Definition}[section]
\newtheorem{thm}[defi]{Theorem}

\newtheorem{corollary}[defi]{Corollary}

\usepackage{latexsym}
\usepackage{amsmath}
\usepackage{amssymb}
\newcommand{\cA}{{\mathcal A}}

\newcommand{\cF}{{\mathcal F}}

\newcommand{\cH}{{\mathcal H}}
\newcommand{\cB}{{\mathcal B}}

\newcommand{\cRH}{{\mathcal {RH}}}

\newcommand{\aA}{{\mathbb A}}

\newcommand{\R}{{\mathbb R}}
\newcommand{\RR}{{\mathbb R}}

\newcommand{\ve}{{\varepsilon}}
\newcommand{\la}{{\langle}}
\newcommand{\ra}{{\rangle}}

\begin{document}

\title[The Microscopic Bidomain Problem]
{On the Microscopic Bidomain Problem\\ with FitzHugh-Nagumo Ionic Transport}

\author{Gieri Simonett}
\address{Department of Mathematics\\
Vanderbilt University\\
Nashville, Tennessee, USA}
\email{gieri.simonett@vanderbilt.edu}

\author{Jan Pr\"uss}
\address{Martin-Luther-Universit\"at Halle-Witten\-berg\\
         Institut f\"ur Mathematik \\
         Theodor-Lieser-Strasse 5\\
         D-06120 Halle, Germany}
\email{jan.pruess@mathematik.uni-halle.de}
\thanks{This work was supported by a grant from the Simons Foundation (\#426729, Gieri Simonett).}



\begin{abstract}
The microscopic  bidomain problem with FitzHhugh-Nagumo ionic transport is studied in the $L_p\!-\!L_q$-framework. Reformulating the problem as a semilinear evolution equation on the interface, local well-posedness is proved in strong as well as in weak settings. We obtain solvability for initial data in the critical spaces of the problem. For dimension $d\leq 3$, by means of energy estimates and a recent result of Serrin type, global existence is shown. Finally, stability of spatially constant equilibria is investigated, to the result that the stability properties of such equilibria  parallel those of the classical FitzHugh-Nagumo system in ODE's. These properties of the bidomain equations are obtained combining recent results on Dirichlet-to-Neumann operators \cite{PrSi16} and \cite{Pru18}, on critical spaces  for parabolic evolution equations \cite{PSW18}, and qualitative theory of evolution equations.
\end{abstract}

\maketitle

\section{Introduction}
 In this paper
we consider the so-called microscopic bidomain system, a well established model 
 describing electrical activity of the heart, 
which is discussed in detail e.g. in the monographs  by Keener and Sneyd ~\cite{KS98} and by Colli Franzone, Pavarino and
Scacchi \cite{CPS14}. To state the problem, consider a bounded domain $\Omega\subset\RR^{d}$, $d\geq2$, with boundary of class $\partial\Omega\in C^{2-}$. This domain consists of two disjoint phases, the intra-cellular space $\Omega_{\sf i}$, and the extra-cellular part $\Omega_{\sf e}$.  The phases are separated by the interface $\Sigma\in C^{2-}$, which is assumed not to intersect the outer boundary $\partial\Omega$, hence $\Sigma=\partial\Omega_{\sf i}$ and $\partial\Omega_{\sf e}=\Sigma \cup\partial\Omega$.
The system is given by
\begin{align}\label{BDE}
\begin{aligned}
\cA(x,D)u &= \phi & \mathrm{in} & \quad \Omega\setminus \Sigma ,  \\
\cB(x,D)u&=\varphi& \mathrm{on} & \quad \partial\Omega,\\
[\![u]\!]=v,\quad -[\![\cB(x,D)u]\!] &=\varphi_\Sigma& \mathrm{on} &\quad \Sigma,\\
\partial_t v - \la\cB(x,D)u\ra &= f(v,w) & \mathrm{on} &\quad \Sigma,\\
\partial_t w &= g(v,w)  & \mathrm{on} & \quad \Sigma,  \\
v(0) = v_0, \quad w(0) &= w_0 &   \mathrm{on} & \quad \Sigma.
\end{aligned}
\end{align}
The function $u$ models the  intra- and extra-cellular electrical potentials, $v =[\![u]\!]= u_{|_{\Omega_{\sf e}}}-u_{|_{\Omega_{\sf i}}} $
is  the transmembrane potential, $w$ is a gating variable, 
 $\nu$ denotes the outward unit normal vector on $\partial \Omega$ and on  $\partial\Omega_{\sf i}$, respectively,
and
 $\la u\ra :=(u_{|_{\Omega_{\sf e}}}+u_{|_{\Omega_{\sf i}}})/2$ is the (algebraic) surface mean of $u$.

Furthermore, assuming sum convention, $\cA(x,D) = - \partial_k{\sf a}^{kl}(x) \partial_l$ denotes a symmetric, uniformly strongly elliptic differential operator in divergence form, and $\cB(x,D) =\nu_k {\sf a}^{kl}(x)\partial_l$ denotes the co-normal derivative operator on $\Sigma\cup \partial\Omega$. We assume for simplicity that the coefficients ${\sf a}^{kl}$ are of class $BU\!C^{1-}(\Omega\setminus\Sigma)$. In particular,  they may jump across the interface $\Sigma$. TThe functions $\phi$, $\varphi$, and $\varphi_\Sigma$ describe exterior sources for the potential $u$ in $\Omega$, which is discontinuous across $\Sigma$.

There are various models describing the ionic transport. Here we choose the most classical one by {\em FitzHugh-Nagumo}, which reads
\begin{equation}
\begin{aligned}
f(v,w)&= -v(v-a)(v-1)-w = -v^3+(a+1)v^2-av-w , \\
g(v,w)&=- bw + cv,
\end{aligned}
\end{equation}
where the constants $a,b$ and $c$ satitsfy
\begin{equation}\label{abc}
0<a<1 \;\text{ and }\; b, c > 0.
\end{equation}

Rigorous mathematical analysis of the above system was pioneered by  Colli-Franzone and Savar\'e~\cite{CS00},
who introduced a variational formulation of the problem, showed global existence and uniqueness of weak  solutions in dimension $d=3$ for the FitzHugh-Nagumo ionic transport, and proved that the solutions are even strong solutions in $L_2$ for initial values $v_0\in H^1_2(\Sigma)$, $w_0\in L_2(\Sigma)$; see also the monograph \cite{CPS14}. These are the only references for the microscopic bidomain problem we are aware of. In particular, an $L_q$-theory of the problem does not seem to exist so far.  

Mathematically, this problem consists of an elliptic transmission problem for the potential $u$, coupled with a dynamic equation on the interface for the jump $v=[\![u]\!]$ and an ODE for the gating variable $w$. This is an interesting, non-standard problem which calls for an approach based on the theory of evolution equations. In fact, we will reformulate the system as a semilinear parabolic evolution equation on the interface $\Sigma$,
employing results for elliptic transmission problems, and in particular results concerning the associated Dirichlet-to-Neumann operator; we refer to Chapter 6.5 of our recent monograph \cite{PrSi16} and to \cite{Pru18}.

To explain this  reformulation, we first consider the elliptic transmission problem
\begin{equation}\label{ETP}
\begin{aligned}
\cA(x,D)u &= \phi & \mathrm{in} & \quad \Omega\setminus \Sigma ,  \\
\cB(x,D)u&=\varphi& \mathrm{on} & \quad \partial\Omega,\\
[\![u]\!]=v,\quad -[\![\cB(x,D)u]\!] &=\varphi_\Sigma& \mathrm{on} &\quad \Sigma,\\
\end{aligned}
\end{equation}
in $L_q(\Omega)$, $1<q<\infty$,  where the data $(\phi,\varphi,\varphi_\Sigma)$ and $v$ are given. Then Theorem 6.5.2 and the subsequent discussion in \cite{PrSi16} shows that
problem \eqref{ETP} admits a solution $u\in H^2_q(\Omega\setminus\Sigma)$, assuming $\phi\in L_q(\Omega)$, $\varphi\in W^{1-1/q}_q(\partial\Omega)$,
$\varphi_\Sigma\in W^{1-1/q}_q(\Sigma)$, and $v\in  W^{2-1/q}_q(\Sigma)$, provided the compatibility condition
$$ \int_\Omega \phi\, dx +\int_{\partial\Omega} \varphi \,d(\partial\Omega)+\int_\Sigma \varphi_\Sigma \,d\Sigma=0$$
is satisfied.
This solution is unique up to a constant and we normalize it by mean 0 over $\Omega$. By linearity, $u=\bar{u}+\tilde{u}$ where $\bar{u}$ denotes the mean zero solution of \eqref{ETP} with $v=0$ and $\tilde{u}=u-\bar{u}$ is only dependent on $v$. The quantity of interest for the dynamics on the interface is
$$ \la\cB(x,D) u\ra  = \la\cB(x,D)\bar{u}\ra  + \cB(x,D)\tilde{u} =: -\Phi + S_{\sf n}v,$$
where $\aA:=-S_{\sf n}$ denotes  the  Dirichlet-to-Neumann operator for the transmission problem \eqref{ETP}, see \cite{PrSi16}, Section 6.5.
Then the problem on the interface $\Sigma$ reads
\begin{equation}\label{dbe}
\begin{aligned}
\partial_t v + \aA v &= \Phi +f(v,w), &&t>0, &&v(0)=v_0,\\
\partial_t w &= g(v,w), &&t>0, &&w(0 )=w_0.
\end{aligned}
\end{equation}
Observe that $v=[\![u]\!] \in W^{2-1/p}_q(\Sigma)$ and $\aA v \in W^{1-1/q}_q(\Sigma)$, hence the latter is the base space for the problem in the strong setting. As $w$ obeys an ode, $w$ can be taken in the same or a better class than $v$, which leads to the following {\em strong formulation} of the problem.
Set $X_0 =W^{1-1/q}_q(\Sigma)^2$, $X_1= W^{2-1/q}_q(\Sigma)\times W^{1-1/q}_q(\Sigma)$ and define $A:X_1\to X_0$ by means of
$$ A = \left[ \begin{array}{cc} \aA +a& 1\\ -c& b\end{array}\right],\quad F({\sf x}) = \left[ \begin{array}{c}- v^3+(a+1)v^2\\0\end{array}\right],$$
where ${\sf x}= [v,w]^{\sf T}$ and ${\sf x}_0= [v_0,w_0]^{\sf T}$ denotes the system variable and its initial value. Then problem \eqref{dbe} can be written as the semilinear parabolic evolution equation in $X_0$
\begin{equation}\label{sleveq}
\partial_t {\sf x} + A{\sf x} = \Phi+F({\sf x}),\quad t>0,\quad {\sf x}(0)={\sf x}_0.
\end{equation}

It should be noted that the classical FitzHugh-Nagumo system appears as a special case of this equation.
To see this,  let $(\phi, \varphi, \varphi_\Sigma)=0$ and $v,w$ be constant on $\Sigma$.

There is also an appropriate weak formulation of problem \eqref{BDE} which reads as follows.
\begin{equation}\label{BDEW}
\begin{aligned}
(a(x)\nabla u|\nabla\Psi)_{L_2(\Omega)} &= \langle \Phi^*|\Psi\rangle, & \mbox{for all } & \Psi\in \dot{H}^1_{q^\prime}(\Omega),\\
[\![u]\!]&=v & \mathrm{on} &\quad \Sigma,\\
\partial_t v +\la\cB(x,D)u\ra  &= f(v,w) & \mathrm{on} &\quad \Sigma,\\
\partial_t w &= g(v,w)  & \mathrm{on} & \quad \Sigma,  \\
v(0) = v_0, \quad w(0) &= w_0 &   \mathrm{on} & \quad \Sigma.
\end{aligned}
\end{equation}
Here $\Phi^*\in \dot{H}^{-1}_q(\Omega\setminus\Sigma):= \dot{H}^{1}_{q^\prime}(\Omega\setminus\Sigma)^*$ and $v_0,w_0$ are given. In this case we have
$u\in H^1_q(\Omega\setminus\Sigma)$ and $\la\cB(x,D)u\ra \in W^{-1/q}_q(\Sigma)$, so as base space we choose $ X_0^{\sf w}= W^{-1/q}_q(\Sigma)\times L_q(\Sigma)$, and the regularity space by standard trace theory becomes $X_1^{\sf w} =  W^{1-1/q}_q(\Sigma)\times L_q(\Sigma)$. In the strong formulation with 
$(\phi, \varphi, \varphi_\Sigma)\in L_q(\Omega)\times L_q(\partial\Omega)\times L_q(\Sigma)$, the functional $\Phi^*$ becomes
$$ \langle \Phi^*|\psi\rangle =\int_\Omega \phi\psi \,dx + \int_{\partial\Omega} \varphi \psi \,d(\partial\Omega)
+\int_\Sigma \varphi_\Sigma  [\![\Psi]\!]  \,d\Sigma,\quad
\psi \in \dot{H}^1_{q^\prime}(\Omega\setminus\Sigma).$$
The weak formulation is as usual derived from the strong formulation by an integration by parts.

In this paper, in order to establish an $L_p\!-\!L_q$-theory, we combine  results on the Dirichlet-to-Neumann operator $S_{\sf n}$ 
derived in \cite{PrSi16}, Section 6.5, and
in \cite{Pru18} with the theory of critical spaces for semi-linear parabolic evolution equations \cite{PSW18},
to obtain local well-posedness of the problem under minimal assumptions on the data in the $L_p\!-\!L_q$-setting with time weights. This gives access to known results on stability of equilibria for parabolic evolution equations and stable and unstable manifolds near hyperbolic equilibria.
A striking consequence of the theory of critical spaces is a very short proof of global existence of strong solutions in the case $d=2$ by using a Serrin type result from \cite{PSW18}. For the cases $d=3$ the proof of global existence is more involved but still simple. It uses energy estimates for the time derivative of the system and the theory of critical spaces.

The plan for this paper is as follows. In Section 2, important results for the underlying linear operator $A$ are derived, which are basic for our approach. Local well-posedness is studied in Section 3 in various functional analytic settings, strong as well as weak, and parabolic regularization is obtained. Section 4 is devoted to elementary energy estimates and global existence in dimensions $d\leq 3$. The final section concerns stability of the homogeneous equilibria for \eqref{sleveq} in the natural state space for the $L_p\!-\!L_q$-theory in case $\Phi\equiv 0$, i.e.\ for $(\phi,\varphi, \varphi_\Sigma)=0$.

\section{The Linear Operator}
From Section 6.5 in \cite{PrSi16} we know that $\aA: W^{2-1/q}_q(\Sigma)\to W^{1-1/q}_q(\Sigma)$ is well-defined and bounded. Theorem 2.1 in \cite{Pru18} shows that it is  pseudo-sectorial in $W^{1-1/q}_q(\Sigma)$ with spectral angle $\phi_\aA=0$. It is symmetric and positive semi-definite in $L_{2}(\Sigma)$, its spectrum consists only of countably many nonnegative eigenvalues of finite algebraic multiplicity, independent of $q$, and  ${\sf N}(\aA)\oplus {\sf R}(\aA)= W^{1-1/q}_q(\Sigma)$.

By symmetry and duality, we see that $\aA: W^{-1/q}_q(\Sigma)\to W^{-1-1/q}_q(\Sigma)$ is also pseudo-sectorial with spectral angle 0 and it has the same properties as $\aA$ has in $W^{1-1/q}_q(\Sigma)$. By means of real interpolation with exponent $q$, this implies that 
$$\aA:W^{s+1}_q(\Sigma) \to W^s_q(\Sigma)$$ 
is again pseudo-sectorial with spectral angle 0 in all spaces $W^s_q(\Sigma)$, for $s\in [-1-1/q,1-1/q]$. But even more is true. Theorem 4.5.4 in \cite{PrSi16} shows that $a+\aA\in \cRH^\infty(W^s_q(\Sigma) )$ with $\cRH^\infty$-angle 0, for all $s\in(-1-1/q,1-1/q)$ and for any $a>0$. As $a+\aA:W^{s+1}_q(\Sigma) \to W^s_q(\Sigma)$ is an isomorphism, these assertions hold also for the limiting values $s=-1-1/q$ and $s=1-1/q$. Note that by regularity, the eigenvalues and eigenspaces of $\aA$ are also independent of $s$. For integer values of $s$, we have to replace $W^m_q(\Sigma)$ by the Besov spaces $B_{qq}^m(\Sigma)$, i.e.\ for $m=0,-1$.

\medskip
Next, we extend $a+\aA$ to $A_0={\rm diag}(a+\aA,b)$ in $X_{0,s}:=W^s_q(\Sigma)\times \cF^s_q(\Sigma)$, where $b>0$. Here $\cF^s_q = W^s_q$ for $s>0$ and $\cF^s_q = L_q$ for $-1<s\leq 0$. Then it is easy to see that $A_0$ is sectorial, invertible, and  $A_0\in \cH^\infty(X_{0,s})$ with $\cH^\infty$-angle $\phi_{A_0}^\infty =0$. Note that the domain of $A_0$ is 
$${\sf D}(A_0)=W^{s+1}_q(\Sigma)\times \cF^s_q(\Sigma)=:X_{1,s}.$$
Now, $A$ as defined before is a bounded perturbation of $A_0$, hence by a perturbation result for the class $\cH^\infty$, we also have $A\in \cH^\infty(X_{0,s})$ and it is sectorial, and then $\phi_A^\infty=\phi_A$, the spectral angle of $A$; see Corollary~3.3.15 in \cite{PrSi16}. To compute the latter, observe that $A-\lambda$ is invertible if and only if $\lambda\neq b$ and $\lambda-a+ c/(\lambda-b)\not\in \sigma(\aA)\subset[0,\infty)$. Actually, $\lambda=b$ belongs to the continuous spectrum of $A$ and is not an eigenvalue; this is an exceptional point. On the other hand, the remaining part of the spectrum $\sigma(A)$ of $A$ consists of eigenvalues of finite algebraic multiplicity, and is determined by the quadratic equation
$$ \lambda^2 -(\mu+a+b) \lambda + (c+b(\mu+a))=0,\quad \mu\in\sigma(\aA).$$
 The solutions of this equation will be either complex, with real part
$${\rm Re}\, \lambda = (\mu+a+b)/2\geq (a+b)/2,$$ or real with
$$\lambda_1>\lambda_2 = \frac{\mu+a+b}{2}\Big( 1- \sqrt{1-4\frac{c+b(\mu+a)}{(\mu+a+b)^2}}\Big)
\geq{\rm min}\Big\{b, \frac{c+ab}{a+b}\Big\}.$$
So we see that
$${\rm Re}\, \sigma(A)\geq \gamma_1:={\rm min}\Big\{b, \frac{a+b}{2},\frac{c+ab}{a+b}\Big\}>0.$$
On the other hand, if there are complex eigenvalues, then their argument $\phi$ satisfies
$$ |{\rm tan}\,\phi|\leq 2\frac{\sqrt{c+b(\mu+a)}}{\mu+a+b} \leq \gamma_3,\quad \mu\geq0,$$
hence $|{\rm arg}\, \sigma(A)|= \phi_A <\pi/2$. This shows that $A$ is sectorial, invertible, and $A\in \cH^\infty(X_{0,s})$ with $\cH^\infty$-angle strictly less than $\pi/2$. In particular, $A$ has maximal $L_p$-regularity on the halfline $\RR_+$, and $-A$ generates an analytic $C_0$-semigroup in $X_{0,s}$ which is exponentially stable. We remind that the most important cases are $s=1-1/q$ for the strong setting and $s=-1/q$ for the weak formulation.
Note that we have
$$ {\sf D}((a+\aA)^\beta) = W^{s+\beta}_q(\Sigma)= {D}_{a+\aA}(\beta,q), \quad {D}_{a+\aA}(\beta,p)= B^{s+\beta}_{qp}(\Sigma),$$
which allows to compute the fractional power space ${\sf D}(A^\beta)$ and the real interpolation spaces ${\sf D}_{A}(\beta,p)$ easily.

\section{Well-Posedness}
In this section we prove local well-posedness of \eqref{sleveq} in the strong as well as in various weak settings.
\subsection{The Strong Approach}
In the strong setting  the base space $X_0$ and the regularity space $X_1$ are given by
$$X_0=W^{1-1/q}_q(\Sigma)^2,\quad X_1= W^{2-1/q}_q(\Sigma)\times W^{1-1/q}_q(\Sigma),$$ 
and
$X_\beta ={\sf D}(A^\beta)= W^{1+\beta-1/q}_q(\Sigma)\times W^{1-1/q}_q(\Sigma)$, $X_{\gamma,\mu}= B^{1+\mu-1/p-1/q}_{q,p}(\Sigma)\times W^{1-1/q}_q(\Sigma).$ 

\medskip
To apply Theorem 1.2 in \cite{PSW18}, we need to estimate the nonlinearity $F$ in an appropriate way. Given $v\in  W^{1+\beta-1/q}_q(\Sigma)$ we let $u\in H^{1+\beta}_q(\Omega)$ denote its Dirichlet extension to $\Omega$; see Appendix B in  \cite[Chapter 10]{PrSi16}.

Then $F(v)$ will be the trace of $F(u)$, hence we need to show that $F:H^{1+\beta}_q(\Omega) \to H^1_q(\Omega)$  is bounded. Then $F:X_\beta\to X_0$ is bounded and polynomial, hence real analytic. In case $q<d$, we estimate as follows.
\begin{align*}
|u^3|_{H^1_q}&\leq c(|u^3|_{L_q} +3|u^2\nabla u|_{L_q})\leq c|u|^2_{L_{2rq}}|u|_{H^1_{r^\prime q}}\leq c |u|_{H^{1+\beta}_q}^3,
\end{align*}
provided $H^{1+\beta}_q(\Omega)\hookrightarrow H^1_{r^\prime q}(\Omega)\hookrightarrow L_{2rq}(\Omega)$ by sharp Sobolev embedding, which requires
$$ 1+\beta -d/q = 1-d/r^\prime q = - d/2rq,\; 1/r+1/r^{\prime}=1,\quad \mbox{i.e.}\; \beta = d/rq= \frac{2}{3}(d/q -1).$$
We have $\beta\in [0,1)$ if $1\le  d/q< 5/2$, which means $2d/5< q\le d.$ In this case the critical weight is
$$\mu_c=1/p+ 3\beta/2 -1/2 = 1/p+ d/q-3/2,$$
see Remark 6.1 in \cite{PSW18}.
Thus we are in the critical case if $d/q>3/2$, i.e.\ $2d/5\leq q< 2d/3$, and the critical spaces are
\begin{equation}\label{critical-strong}
X_{crit} = X_{\gamma,\mu_c}= B^{(d-1)/q-1/2}_{qp}(\Sigma)\times W^{1-1/q}_q(\Omega).
\end{equation}
 For $q> d$ we have $H^1_q(\Omega)\hookrightarrow L_\infty(\Omega)$, in this case
$F: H^1_q(\Omega)\to H^1_q(\Omega)$ is bounded, hence $\beta=0$ for $q>d$. In the borderline case $q=d$, we may choose $\beta>0$ as small as we want.
This shows that the well-posedness range is $q> 2d/5$.
Therefore, Theorem 1.2 in \cite{PSW18} yields the following results on strong local well-posedness of \eqref{sleveq}.

\begin{thm}\label{str-set}
Let $d\geq2$,  $p,q\in (1,\infty)$ with $1/p+d/q\leq 5/2$, 
 let $\mu\in (1/p,1]$ satisfy
$$ \mu\geq \mu_c = \frac{1}{p} + \frac{d}{q} -\frac{3}{2},$$ 
and let $\Phi\in L_{p,\mu,loc}([0,\infty); W^{1-1/q}_q(\Sigma))$.
Then for each initial value  ${\sf x}_0=(v_0,w_0)$ with regularity
$$(v_0,w_0)\in B^{1+\mu-1/p-1/q}_{qp}(\Sigma)\times W^{1-1/q}_q(\Sigma),$$
there exists $\tau>0$ and a unique solution $(v,w)$ of \eqref{sleveq} in the class
\begin{equation*}
\begin{split}
&v\in H^1_{p,\mu}((0,\tau); W^{1-1/q}_q(\Sigma))\cap L_{p,\mu}((0,\tau); 
 W^{2-1/q}_q(\Sigma))\hookrightarrow C([0,\tau]; B^{1+\mu-1/p-1/q}_{qp}(\Sigma)), \\
& w\in H^2_{p,\mu}((0,\tau);W^{1-1/q}_q(\Sigma))\hookrightarrow C^{1+\mu-1/p}([0,\tau];W^{1-1/q}_q(\Sigma)).
\end{split}
\end{equation*}
The solution depends continuously on the data in the corresponding spaces, and it exists on a maximal time interval $[0,t_+({\sf x}_0))$, which is characterized by
$$t_+<\infty\; \Leftrightarrow \; \lim_{t\to t_+} v(t)  \mbox{ does not exist in }  B^{1+\mu-1/p-1/q}_{qp}(\Sigma).$$
In addition, $v$ regularizes instantly if $1/p+d/q< 5/2$, it belongs to
$$ v\in H^1_{p,loc}((0,t_+);W^{1-1/q}_q(\Sigma))\cap L_{p,loc}((0,t_+); W^{2-1/q}_{q}(\Sigma))
\hookrightarrow C((0,t_+); B^{2-1/p-1/q}_{qp}(\Sigma)) ,$$
and we have
$ u\in L_{p,loc}((0, t_+); H^2_q(\Omega\setminus\Sigma)).$
\end{thm}
\noindent
The improved time regularity of $w$ is obtained by using directly the equation for $w$. Note that $w$ does not gain spatial regularity, as there is no diffusion in the $w$-equation.
We observe that in the critical case $2d/5< q<2d/3$ we may choose $\mu=\mu_c$, which leads to well-posedness in the corresponding critical spaces.

\subsection{The Weak Appoach}
In the weak setting we choose
$$X_0^{\sf w}=W^{-1/q}_q(\Sigma)\times L_q(\Sigma),\quad X_1^{\sf w}=W^{1-1/q}_q(\Sigma)\times L_q(\Sigma),$$ 
which yields
$X_\beta^{\sf w} = W^{\beta-1/q}_q(\Sigma)\times L_q(\Sigma)$ and $X_{\gamma,\mu}^{\sf w}= B_{qp}^{\mu-1/p-1/q}(\Sigma)\times L_q(\Sigma)$. This choice leads to $u\in H^1_q(\Omega\setminus\Sigma)$. To find the optimal value of $\beta^{\sf w}\in[0,1)$ we consider the cubic part $F_0: H_q^{\beta-s}(\Sigma)^3\to H^{-s}_q(\Sigma)$. We estimate  in the following way. Fix $\phi\in H^s_{q^\prime}(\Sigma)$; then
$$ |\langle F_0(v)|\phi\rangle|\leq |v|_{L_{3r}}|\phi|_{L_{r^\prime}},$$
with optimal choices by sharp Sobolev embedding.
Hence,
$$ s+\frac{d-1}{q}= \frac {d-1}{r}<d-1,\quad \beta =\frac{2}{3}(s+\frac{d-1}{q}),$$
which is feasible provided $r>1$ and $\beta -s =(d-1) (1/q- 1/3r)\geq 0$. These conditions are equivalent to $s<(d-1)(1-1/q)$ and $s<2(d-1)/q$, while
$\beta<1$ means $s<3/2-(d-1)/q$. By means of real interpolation between $s=1/q-\ve$ and $s=1/q+\ve$, this yields $F: X^{\sf w}_{\beta^{\sf w}}\to X_0^{\sf w}$ polynomial and bounded, hence real analytic, provided
$$ \frac{1}{q} < \min\Big\{ 2\frac{d-1}{q}, \frac{d-1}{q^\prime}, \frac{3}{2}-\frac{d-1}{q}\Big\},\quad \beta^{\sf w} = \frac{2d}{3q},$$
which are equivalent to $ q>2 $ for $d=2$ and $q> 2d/3$ for $d\geq 3$. The critical weight becomes
$$ \mu_c^{\sf w} =\frac{1}{p} + \frac{3\beta^{\sf w}}{2}-\frac{1}{2} = \frac{1}{p} +\frac{d}{q} - \frac{1}{2},$$
hence we are in the critical range if $ 2\max\{1,d/3\}<q<2d$. Applying Theorem 1.2 in \cite{PSW18} this yields

\begin{thm}\label{weak-set}
Let $d\geq2$,   $p,q\in (1,\infty)$ with $1/p+d/q\leq 3/2$,  
and let $\mu\in (1/p,1]$ satisfy $q>2$ in case $d=2$ and
$$ \mu\geq \mu_c^{\sf w} = \frac{1}{p} + \frac{d}{q} -\frac{1}{2}$$
and let $\Phi\in L_{p,\mu,loc}([0,\infty); W^{-1/q}_q(\Sigma)).$
Then for each initial value ${\sf x}_0=(v_0,w_0)$ with regularity
$(v_0,w_0) \in B^{\mu-1/p-1/q}_{qp}(\Sigma)\times L_q(\Sigma) $
there exists $\tau>0$ and a unique solution $(v,w)$ of \eqref{sleveq} in the class
\begin{equation*}
\begin{split}
&v\in H^1_{p,\mu}((0,\tau); W^{-1/q}_q(\Sigma))\cap L_{p,\mu}((0,\tau); W^{1-1/q}_q(\Sigma))
\hookrightarrow C([0,\tau]; B^{\mu-1/p-1/q}_{qp}(\Sigma)), \\
& w\in  W^{2-1/q}_{p,\mu}((0,\tau);L_q(\Sigma))\hookrightarrow C^{1+\mu-1/p-1/q}([0,\tau];L_q(\Sigma)).
\end{split}
\end{equation*}
The solution depends continuously on the data in the corresponding spaces, and it exists on a maximal time interval $[0,t_+^{\sf w}({\sf x}_0))$, which is characterized by
$$t_+^{\sf w}<\infty\; \Leftrightarrow \; \lim_{t\to t_+^{\sf w}} v(t) \mbox{ does not exist in }  B^{\mu-1/p-1/q}_{qp}(\Sigma).$$
In addition, $v$ regularizes instantly if $1/p+d/q< 3/2$; it belongs to
$$ v\in H^1_{p,loc}((0,t_+^{\sf w});W^{-1/q}_q(\Sigma))\cap L_{p,loc}((0,t_+^{\sf w}); W^{1-1/q}_{q}(\Sigma))\hookrightarrow C((0,t_+^{\sf w}); B^{1-1/p-1/q}_{qp}(\Sigma)),$$
and we have
$ u\in L_{p,loc}((0, t_+^{\sf w}); H^1_q(\Omega\setminus\Sigma)).$
\end{thm}
\noindent
We observe that in the critical case $2<q<4$ for $d=2$, and $2d/3<q<2d$ for $d\geq3$, we may choose $\mu=\mu_c$, which leads to well-posedness in the critical spaces 
$$X_{crit}^{\sf w} = B_{qp}^{(d-1)/q-1/2}(\Sigma)\times L_q(\Sigma).$$ 

We add one small remark on the characterization of the maximal time of existence $t_+^{\sf w}$. If  $\mu-1/p-1/q>0$, then $B^{\mu-1/p-1/q}_{qp}(\Sigma)\hookrightarrow L_q(\Sigma)$, so the result follows from Corollary~2.3 in \cite{PSW18}. For $\mu\leq 1/p+1/q$, we first obtain from the equation for $w$ that $\lim_{t\to t_+^{\sf w}} w(t)$ exists in  $B^{\mu-1/p-1/q}_{qp}(\Sigma)$, hence in $W^{-1/q}_q(\Sigma)$. Then looking at the equation for $v$ we obtain from the Corollary just cited $v\in L_{p,\mu}((0,t_+^{\sf w});W^{1-1/q}_q(\Sigma))$, and so $\lim_{t\to t_+^{\sf w}} w(t)$ exists also in $L_q(\Sigma)$.

\subsection{An Extended Weak Setting}
To cover the gaps and extend the range of criticality in the strong and weak settings, in particular to cover the case $q=2$ for dimensions $d=2,3$, we introduce a general weak setting. For this purpose, we choose $\kappa\geq0$, and set 
$$X_0^\kappa=B^{-\kappa/q}_{qq}(\Sigma)\times L_q(\Sigma),
\quad  X_1^\kappa=B^{1-\kappa/q}_{qq}(\Sigma)\times L_q(\Sigma),$$ 
which yields
$X_\beta^\kappa = W^{\beta-\kappa/q}_q(\Sigma)\times L_q(\Sigma)$ and 
$X_{\gamma,\mu}^\kappa= B_{qp}^{\mu-1/p-\kappa/q}(\Sigma)\times L_q(\Sigma).$ 
This choice leads to $u\in H^{1-(\kappa-1)/q}_q(\Omega\setminus\Sigma)$, provided $\kappa<q$.

The nonlinearity is estimated in the same way as in the previous subsection, to the result that
$$ \beta^\kappa = 2(d+\kappa-1)/3q,\quad \mu_c^\kappa = 1/p+ (d+\kappa-1)/q -1/2,$$
which is feasible if $\beta^\kappa>\kappa/q$, equivalently $\kappa< 2(d-1)$, with the constraints
\begin{equation}\label{constraints-kappa}
1+\kappa/(d-1)<q,\quad 2(d+\kappa-1)/3<q,\quad \mbox{and for criticality} \; q< 2(d+\kappa-1).
\end{equation}
Observe that these conditions imply $0<\kappa<q$.
 To see this let us assume, to the contrary, that $q<\kappa$. Then the second condition in 
\eqref{constraints-kappa} yields $2(d-1)<q$.
But then $2(d-1)<q<\kappa$, contradicting the constraint $\kappa<2(d-1).$ 

Applying once more Theorem 1.2 in \cite{PSW18}, this yields the following result.
\begin{thm}\label{ext-weak-set}
Let  $d\geq2$,  $\kappa\in [0,2(d-1))$, $p,q\in (1,\infty)$ with
$1/p+(d+\kappa-1)/q\leq 3/2,$
  and let $\mu\in (1/p,1]$ satisfy
$1+\kappa<q$ for $d=2$, and
$$ \mu\geq \mu_c^\kappa = \frac{1}{p} + \frac{d+\kappa-1}{q} -\frac{1}{2}$$
and let $\Phi\in L_{p,\mu,loc}([0,\infty); B^{-\kappa/q}_{qq}(\Sigma))$.
Then for each initial value ${\sf x}_0=(v_0,w_0)$ with regularity 
$(v_0,w_0) \in B^{\mu-1/p-\kappa/q}_{qp}(\Sigma)\times L_q(\Sigma),$
there exists $\tau>0$ and a unique solution $(v,w)$ of \eqref{sleveq} in the class
\begin{equation*}
\begin{split}
&v\in H^1_{p,\mu}((0,\tau); B^{-\kappa/q}_{qq}(\Sigma))\cap L_{p,\mu}((0,\tau); B^{1-\kappa/q}_{qq}(\Sigma))\hookrightarrow C([0,\tau]; B_{qp}^{\mu-1/p-\kappa/q}(\Sigma)), \\
& w\in W^{2-\kappa/q}_{p,\mu}((0,\tau);L_q(\Sigma))\hookrightarrow C^{1+\mu-1/p-\kappa/q}([0,\tau];L_q(\Sigma)).
\end{split}
\end{equation*}
The solution depends continuously on the data in the corresponding spaces, and it exists on a maximal time interval $[0,t_+^\kappa({\sf x}_0))$, which is characterized by
$$t_+^\kappa<\infty \; \Leftrightarrow\; \lim_{t\to t_+^\kappa} v(t) \mbox{ does not exist in }  B^{\mu-1/p-\kappa/q}_{qp}(\Sigma).$$
In addition, $v$ regularizes instantly; it belongs to
$$ v\in H^1_{p,loc}((0,t_+^\kappa);B^{-\kappa/q}_{qq}(\Sigma))\cap L_{p,loc}((0,t_+^\kappa); B^{1-\kappa/q}_{qq}(\Sigma))\hookrightarrow C((0,t_+^\kappa); B^{1-1/p-\kappa/q}_{qp}(\Sigma)),$$
and we have
$ u\in L_{p,loc}((0, t_+^\kappa); H^{1-(\kappa-1)/q}_q(\Omega\setminus\Sigma)).$
\end{thm}
\noindent
{\bf Remarks.} {\bf (i)} The full range for $q$  in the general weak setting is $1<q<6$ for $d=2$ and $2(d-1)/3<q<6(d-1)$, choosing $\kappa$ appropriately, and the critical spaces are given by  
\begin{equation}\label{crtical-spaces-kappa}
X_{crit}= B_{qp}^{(d-1)/q-1/2}(\Sigma)\times L_q(\Sigma).
\end{equation} 
Observe that these spaces are increasing with $p$. The Sobolev exponents for the $v$-part are $-1/2$, independent of $p,q$. This implies that the `largest critical space' arising from $L_p\!-\!L_q$ theory for the $v$-part is $B_{6(d-1),\infty}^{-1/3}(\Sigma)\times L_q(\Sigma)$. \vspace{2mm}\\
 {\bf (ii)} For $p=q=2$ and $\kappa=0$ we have $\mu_c^0= \frac{d-1}{2}$, $\beta^0= (d-1)/3$, and the spaces $X_j^\kappa$ become
$$ X_0^0=L_2(\Sigma)^2,\; X_1^0= H^1_2(\Sigma)\times L_2(\Sigma),\; X_\beta^0 
=H^{(d-1)/3}_2(\Sigma)\times L_2(\Sigma),\;  X_{\gamma,\mu}^0= H^{\mu-1/2}_2(\Sigma)\times L_2(\Sigma).$$
In particular, $u\in H^{3/2}_2(\Omega\setminus\Sigma)$, and $X_{crit} = H^{(d-2)/2}_2(\Sigma)\times L_2(\Sigma)$ for dimensions 2 and~3. 
Here we note that for $d=2$, the constraint for criticality $q<2(d+\kappa-1)$ in \eqref{constraints-kappa}
is not met. However, any $\mu\in (\mu^\kappa_c,1]=(1/2,1]$ is admissible.
We summarize this situation by saying that $d=2$ is borderline subcritical.
If $d=3$, then all constraints listed in \eqref{constraints-kappa} and in Theorem~\ref{ext-weak-set} are satisfied
with $\mu=\mu^\kappa_c$. Hence, for $d=3$ we are in the critical case.
\vspace{2mm}\\
{\bf (iii)} For $p=q=2$ and $\kappa=1$ we have $\mu_c^{1} = d/2$, $\beta =d/3$, 
$$ X_0^1=H^{-1/2}_2(\Sigma)\times H^{1/2}_2(\Sigma),\; X_1^1= L_2(\Sigma)^2,\;
X^1_\beta = H^{d/3-1/2}_2(\Sigma)\times L_2(\Sigma), $$
 and $X^1_{\gamma,\mu}= H^{\mu-1}_2(\Sigma)\times L_2(\Sigma).$
In particular $u\in H^1_2(\Omega\setminus\Sigma)$ and $X_{crit}= H^{(d-2)/2}_2(\Sigma)\times L_2(\Sigma)$. 
Here the constraints in \eqref{constraints-kappa} do not hold for $d=2,3$.
Indeed, one notes that the first constraint in \eqref{constraints-kappa} fails for $d=2$, 
while the second (which reflects the condition $\beta<1) $ and the third (which reflects the condition $\mu\le 1)$  both fail for $d=3$.

\subsection{Further Regularity}
Parabolic regularization implies that the regularity of the solution is actually better than stated in the previous results, provided $\Phi$ has also more regularity. Here we assume $\Phi\in L_{p,\mu,loc}([0,\infty);W^{1-1/q}_q(\Sigma))$, as well as $w_0\in W^{1-1/q}_q(\Sigma)$. 
Let 
$$0<\kappa<2(d-1),\quad \mu_c=1/p+(d+\kappa-1)/q-1/2\in (1/p,1],$$ 
and suppose that $v_0$ belongs to the critical space $ B_{qp}^{(d-1)/q-1/2}(\Sigma)$.  Then for any $\delta_1>0$ t
he solution satisfies ${\sf x}(\delta_1)=(v(\delta_1),w(\delta_1))\in B^{1-1/p-\kappa/q}_{qp}(\Sigma)\times L_q(\Sigma)$. As
$$1-1/p -\kappa/q> \mu_1 -1/p -\ve/q,\quad \mbox{i.e.} \quad 1/p<\mu_1<1-(\kappa-\ve)/q,$$
for $\ve>0$ sufficiently small and appropriate $\mu_1\in (1/p,1]$, 
we obtain
$$v(\delta_1)\in B^{\mu_1-1/p-\ve/q}_{qp}= (W^{-\ve/q}_q(\Sigma),W^{1-\ve/q}_q(\Sigma))_{\mu_1-1/p,p}.$$
Hence, the solution can be continued in $W^{-\ve/q}_q(\Sigma)\times L_q(\Sigma)$. Then, after any second small time $\delta_2>0$, we have $ v(\delta_1+\delta_2)\in B^{1-1/p-\ve/q}_{qp}(\Sigma)$. 
As this space embeds into 
$$B^{1+\mu_2-1/p-1/q}_{qp}(\Sigma)=(W^{1-1/q}_q,W^{2-1/q}_q)_{\mu_2-1/p,p},$$
provided $1/p<\mu_2<(1-\ve)/q$, this yields by the equation for $w$
$${\sf x}( \delta_1+\delta_2)\in B^{1+\mu_2-1/p-1/q}_{qp}(\Sigma)\times W^{1-1/q}_q(\Sigma),$$
and therefore the solution can be continued as a strong solution in $X_0=W^{1-1/q}_q(\Sigma)^2$. We emphasize that, regarding the critical spaces, $p$ is a kind of `play parameter'. We may choose it as large as we want, as $B^s_{qp}(\Sigma)\subset B^s_{qp_1}(\Sigma)$, for all $s$ and $q$, provided $p_1\geq p$. These considerations yield the following result.

\begin{thm}\label{regularity}
Let $d\geq2$, $p,q\in (1,\infty)$ 
with $1<q<6$ for $d=2$, $2(d-1)/3<q<6(d-1)$ for $d\geq3$ such that
$$ \frac{1}{p}<\frac{1}{q},\quad \frac{1}{p}+\frac{d-1}{q}<3/2,$$
and let $\Phi\in L_{p,loc}([0,\infty)); W^{1-1/q}_q(\Sigma))$.
Then for each initial value $(v_0,w_0)$ with regularity 
$$(v_0,w_0) \in B^{(d-1)/q-1/2}_{qp}(\Sigma)\times W^{1-1/q}_q(\Sigma),$$
the unique  solution $(v,w)$ of \eqref{sleveq} belongs to the class
\begin{equation*}
\begin{split}
&v\in H^1_{p,loc}((0,t_+); W^{1-1/q}_q(\Sigma))\cap L_{p,loc}((0,t_+); W^{2-1/q}_q(\Sigma))
\hookrightarrow C((0,t_+); B^{2-1/p-1/q}_{qp}(\Sigma)),\\
& w\in H^{2}_{p,loc}((0,t_+);W^{1-1/q}_q(\Sigma))\hookrightarrow C^{2-1/p}([0,\tau];W^{1-1/q}_q(\Sigma)),
\end{split}
\end{equation*}
and 
$ u \in L_{p,loc}((0,t_+); H^2_q(\Omega\setminus\Sigma)).$
\end{thm}
\section{Energy Estimates}
In this section we prove global existence of the solutions in dimensions $d=2$ and $d=3$.
In the following, we fix $\tau\in (0,t_+)$, where $t_+$ is the maximal time of existence of a solution.

\subsection{The Energy Estimate}
Taking the inner product in $L_2(\Sigma)$ of the equation for $v$ with the solution $v$ we obtain
\begin{align*}
(\partial_tv|v)_{L_2(\Sigma)} + \langle(a+\aA)v|v\rangle + |v|_{L_4(\Sigma)}^4 +(w|v)_{L_2(\Sigma)} = \langle\Phi|v\rangle +(a+1)(v^2|v)_{L_2(\Sigma)},
\end{align*}
and from the equation for $w$ we obtain
$$ (\partial_tw|w)_{L_2(\Sigma)} = - b|w|_{L_2(\Sigma)}^2 +c(w|v)_{L_2(\Sigma)}.$$
Combining these identities yields with
$\psi = \frac{1}{2}|v|^2_{L_2(\Sigma)} + \frac{1}{2c} |w|^2_{L_2(\Sigma)}$
the {\em energy balance}
\begin{align}\label{energy}
 \partial_t \psi  + |(a+\aA)^{1/2}v|^2_{L_2(\Sigma)} + |v|^4_{L_4(\Sigma)} + \frac{b}{c}|w|^2_{L_2(\Sigma)}
= \langle\Phi|v\rangle +(a+1) (v^2|v)_{L_2(\Sigma)}.
\end{align}
As ${\sf D}((a+\aA)^{1/2})= H^{1/2}_2(\Sigma)$ in $L_2(\Sigma)$, there is a constant $c_0>0$ such that
$$ c_0|v|_{H^{1/2}_2(\Sigma)} \leq |(a+\aA)^{1/2}v|^2_{L_2(\Sigma)},$$
and the right hand side of \eqref{energy} is bounded by
\begin{align*}
&|\Phi|_{H^{-1/2}_2(\Sigma)}|u|_{H^{1/2}_2\Sigma)}+ \frac{1}{2}|v|_{L_4(\Sigma)}^4 + \frac{(a+1)^2}{2}|v|_{L_2(\Sigma)}^2\\
&\leq \frac{c_0}{2}|v|_{H^{1/2}_2( \Sigma)}^2 + \frac{1}{2}|v|_{L_4(\Sigma)}^4  + (a+1)^2 \psi + \frac{1}{2c_0}|\Phi|_{H^{-1/2}_2(\Sigma)}^2.
\end{align*}
Therefore, we have
$$\partial_t \psi\leq  (a+1)^2 \psi + \frac{1}{2c_0}|\Phi|_{H^{-1/2}_2(\Sigma)}^2,$$
hence Gronwall's inequality implies
$$ \psi(t) \leq e^{(a+1)^2t}[\psi(0) + (2c_0)^{-1}\int_0^t |\Phi(s)|^2_{H^{-1/2}_2(\Omega\setminus\Sigma)}\,ds],\quad t\in (0, \tau).$$
In total, this yields the energy bounds
\begin{align}\label{enbds}
 v,w,\dot{w}\in L_\infty((0,\tau);L_2(\Sigma)),\quad v\in L_2((0,\tau);H^{1/2}_2(\Sigma)\cap L_4(\Sigma)),
\end{align}
for all space dimensions $d\geq2$, provided
$$ v_0,w_0\in L_2(\Sigma),\; \Phi\in L_2((0,\tau);H^{-1/2}_2(\Sigma)).$$
Now suppose that $d=2$ and
$$ (v_0,w_0)\in B_{qp}^{1/q-1/2}(\Sigma)\times L_q(\Sigma),\quad \Phi\in L_{p,loc}([0,\infty);W^{-1/q}_q(\Sigma)),$$
for $p\geq2$ and $2< q<6$.
Then the unique solution ${\sf x}=[v,w]^{\sf T}$ according to Theorem 3.3 and parabolic regularization with $\kappa=2-\ve$ for $3-\ve<q<6-2\ve$ resp.\ Theorem 3.2 for $2<q<4$ satisfies
$$ v(\delta)\in B^{1-1/p-1/q}_{qp}(\Sigma)\hookrightarrow L_2(\Sigma),$$
for any $\delta>0$, and with $W^{-1/q}_q(\Sigma)\hookrightarrow H^{-1/2}_2(\Sigma)$ we also have 
$\Phi\in L_{2,loc}([\delta,\infty);H^{-1/2}_2(\Sigma))$. 
\cite[Proposition A.1]{PSW18}, 
Sobolev embedding, and  the embedding $H^s_q(\Sigma)\hookrightarrow W^s_q(\Sigma)$, valid for $q\ge 2$, implies
$$ L_{\infty}(J; L_2(\Sigma))\cap L_2(J; H^{1/2}_2(\Sigma))
\hookrightarrow L_p(J; H^{1/p}_2(\Sigma))\hookrightarrow L_p(J;W_q^{1/p+1/q-1/2}(\Sigma))
\quad p,q\ge 2,$$
for $J\subset \R$. 
Therefore, the energy estimate yields
$$v\in L_p([\delta,\tau); W_q^{1/p+2/q-1/2-1/q}(\Sigma))=L_p([\delta,\tau); W_q^{\mu_c^{\sf w}-1/q}(\Sigma)),$$
and the equation for $w$ yields $w\in L_p([\delta,\tau); L_q(\Sigma))$. This shows ${\sf x}\in L_p([\delta;\tau)X^{\sf w}_{\mu_c^{\sf w}})$, 
where $X^{\sf w}_{\mu_c^{\sf w}}$ denotes the complex interpolation space $(X^{\sf w}_0,X^{\sf w}_1)_{\mu_c^{\sf w}}$. 
Since $\tau\in (0,t_+)$ can be chosen arbitrarily, the Serrin type criterion Theorem 2.4 in \cite{PSW18} yields global existence. This way, we have shown the following result on global existence in the case $d=2$.

\begin{thm}
Let $d=2$, $p\geq 2$, $2< q<6$, $\Phi\in L_{p,loc}([0,\infty); W^{-1/q}_q(\Sigma))$.
Then for each initial value ${\sf x}_0=(v_0,w_0)$ with regularity
$$ (v_0,w_0) \in B^{1/q-1/2}_{qp}(\Sigma)\times L_q(\Sigma),$$
the unique solution $(v,w)$ of \eqref{sleveq} exists globally in the class
\begin{equation*}
\begin{split}
&v\in H^1_{p,loc}((0,\infty); W^{-1/q}_q(\Sigma))\cap L_{p,loc}((0,\infty); W^{1-1/q}_q(\Sigma))
\hookrightarrow C((0,\infty); B^{1-1/p-1/q}_{qp}(\Sigma)),\\
& w\in W^{2-1/q}_{p,loc}((0,\infty);L_q(\Sigma))\hookrightarrow C^{2-1/p-1/q}((0,\infty);L_q(\Sigma)),
\end{split}
\end{equation*}
and  $u\in L_{p,loc}((0,\infty); H^1_q(\Omega\setminus\Sigma)).$
\end{thm}

The case $q=2$ needs further considerations, as Theorem 3.2 does not apply. Choosing $\kappa\in (0,1)$, 
Theorem 3.3 yields a unique strong solution with
$$v\in H^1_{p,loc}((0,t_+); H^{-\kappa/2}_2(\Sigma))\cap L_{p,loc}((0,t_+); H^{1-\kappa/2}_2(\Sigma))\hookrightarrow C((0,t_+); B^{1-1/p-\kappa/2}_{2p}(\Sigma)),$$
provided
$$ v_0\in B_{2p}^{0}(\Sigma),\; w_0\in L_2(\Sigma),\quad \Phi\in L_{p,loc}([0,\infty);H^{-\kappa/2}_2(\Sigma)).$$
As 
$B^{1-1/p-\kappa/2}_{2p}(\Sigma))\hookrightarrow L_2(\Sigma)$, the energy bounds are available on $[\delta,\tau)$ 
for any $\tau\in (\delta, t_+)$.
Hence we may conclude as before that the solution is in fact global. These considerations yield
\begin{corollary} In the situation of Theorem 4.1, let $d=q=2$, $p\geq2$, $\kappa\in [0,1)$, and asssume
$$ (v_0,w_0)\in B_{2p}^{0}(\Sigma)\times L_2(\Sigma),\;\Phi\in L_{p,loc}([0,\infty);H^{-\kappa/2}_2(\Sigma)).$$
Then the unique solution $(v,w)$ of \eqref{sleveq} is global and has regularity
\begin{equation*}
\begin{split}
&v\in H^1_{p,loc}((0,\infty); H^{-\kappa/2}_2(\Sigma))\cap L_{p,loc}((0,\infty); H^{1-\kappa/2}_2(\Sigma))
\hookrightarrow C((0,\infty); B^{1-1/p-\kappa/2}_{2p}(\Sigma)),\\
& w\in  W^{2-1/q}_{p,loc}((0,\infty);L_2(\Sigma))\hookrightarrow C^{2-1/p-\kappa/2}((0,\infty);L_2(\Sigma)),
\end{split}
\end{equation*}
and $u\in L_{p,loc}((0,\infty); H^{1+(1-\kappa)/2}_2(\Omega\setminus\Sigma)).$
\end{corollary}

\subsection{Higher Energy Estimate}
To prove the second energy estimate, we take the inner product of the equation for $v$ with $\dot{v}$ in $L_2(\Sigma)$. This yields
\begin{align*}
|\dot{v}|^2_{L_2} + \frac{d}{dt}\big[\frac{1}{2} |(a+\aA)^{1/2}v|^2_{L_2}+ \frac{1}{4}|v|_{L_4}^4 -\frac{a+1}{3}(v^2|v)_{L_2}\big]= (w|\dot{v})_{L_2}+\langle \Phi|\dot{v}\rangle.
\end{align*}
Integrating from 0 to $t$ this yields
\begin{align*}
|\dot{v}|_{L_2((0,t);L_2)}^2 +  \frac{c_0}{2}|v(t)|_{H^{1/2}_2}^2 &+  \frac{1}{4}|v(t)|_{L_4}^4\leq |\int_0^t\langle\Phi(s)\dot{v}(s)\rangle\, ds| + |w|_{L_2((0,t_+);L_2)}|\dot{v}|_{L_2((0,t);L_2)}\\ 
&+ C|v_0|_{H^{1/2}_2}^2 + \frac{1}{4}|v_0|_{L_4}^4 + \frac{a+1}{3}\big[ |v(t)|_{L_4}^2|v(t)|_{L_2}+ |v_0|_{L_4}^2|v_0|_{L_2}^2\big].
\end{align*}
Next we integrate by parts
$$ \int_0^t  \langle\Phi(s)|\dot{v}(s)\rangle\,ds= -\int_0^t\langle\dot{\Phi}(s)|v(s)\rangle \,ds+ \langle\Phi(t)|v(t)\rangle -\langle\Phi(0)|v_0\rangle,$$
and estimate
$$\!\int_0^t \!\! \langle\Phi(s)|\dot{v}(s)\rangle ds|\leq |\dot{\Phi}|_{L_2((0,t);H^{-1/2}_2)}|v|_{L_2((0,t);L_2)}  +  |\Phi(t)|_{H^{-1/2}_2}|v(t)|_{H^{\frac{1}{2}}_2} +  |\Phi(0)|_{H^{-1/2}_2}|v_0|_{H^{1/2}_2}.$$
Employing Young's inequality, and using the elementary energy estimate from the previous subsection, this yields
$$ v\in L_\infty((0,\tau);H^{1/2}_2(\Sigma)\cap L_4(\Sigma)),\quad v,w,\dot{w}\in H^1_2((0,\tau);L_2(\Sigma)),$$
in all dimensions, provided
$$v_0\in H^{1/2}_2(\Sigma)\cap L_4(\Sigma),\quad w_0\in L_2(\Sigma),\quad \Phi\in H^1_{2,loc}([0,t_+); H^{-1/2}_2(\Sigma)).$$
To derive the third energy estimate, we choose a non-decreasing cut-off function $\chi\in C^\infty(\RR)$ with $\chi(t)=0$ for $t\leq\delta/3$, $\chi(t)=1$ for $t\geq 2\delta/3$,
differentiate the equations w.r.t.\ time $t$, and set $\bar{v}=\chi\dot{v}$, $\bar{w}=\chi\dot{w}$. Then $(\bar{v},\bar{w})$ satisfies
\begin{align*}
 \partial_t\bar{v} +(a+\aA)\bar{v} +3v^2 \bar{v}+ \bar{w}&= \chi\dot{\Phi} +2(a+1) v\bar{v} +\dot{\chi}\dot{v}\\
 \partial_t\bar{w} + b\bar{w}-c\bar{w} &= \dot{\chi}\dot{w}.
 \end{align*}
Applying the energy estimate as in  the previous subsection and using $\dot{v},\dot{w}\in L_2((0,\tau);L_2(\Sigma))$ yields the bounds
\begin{align}\label{enbds2}
 v,w\in W^1_\infty((\delta,\tau);L_2(\Sigma)),\quad v\in H^1_2((\delta, \tau); H^{1/2}_2(\Sigma)),\quad v\dot{v} \in L_2((\delta,\tau);L_2(\Sigma)),
\end{align}
provided
$$ v_0\in H^{1/2}_2(\Sigma)\cap L_4(\Sigma),\; w_0\in L_2(\Sigma),\; \Phi\in H^1_{2,loc}([0,t_+);H^{-1/2}_2(\Sigma)).$$
This estimate is also valid in all dimensions $d\geq2$.

By assumption, $\Phi\in H^1_{p,loc}([0,\infty); W^{-1/q}_q(\Sigma))$ and hence $\Phi\in C([0,\infty); W^{-1/q}_q(\Sigma))$
by embedding. Therefore, we may use Theorem 3.2 with $r\geq p$ sufficiently large instead of $p$, to obtain a local solution in the right class as soon as $2<q<6$. Then by regularity of $\Phi$ we may apply Theorem 5.2.1 in \cite{PrSi16}, to obtain
$$ v\in H^2_{p,loc} ([\delta,t_+^{\sf w}); W^{-1/q}_q(\Sigma))\cap  H^1_{p,loc} ([\delta,t_+^{\sf w}); W^{1-1/q}_q(\Sigma)),$$
for any $\delta>0$, hence $v(\delta)\in  W^{1-1/q}_q(\Sigma)\hookrightarrow L_2(\Sigma)$. In the range $6-\ve<q<12-2\ve$ we employ first Theorem 3.3 with
$$v(\delta_1)\in B^{\mu_1- 1/p-1/q}_{qp}(\Sigma)=(W^{-1/q}_q(\Sigma), W^{1-1/q}_q(\Sigma))_{\mu_1-1/p,p}, $$ 
for some $\mu_1\in (1/2,1]$.
Hence we may continue the solution in $W^{-1/q}_q(\Sigma)$ and may conclude as before that $v(\delta_1+\delta_2)\in L_2(\Sigma)$. We again observe that
$W^{-1/q}_q(\Sigma)\hookrightarrow H^{-1/2}_2(\Sigma)$. Therefore, these solutions satisfy the energy estimates \eqref{enbds2}.

For $d=3$ and $p,q\geq 2$ we have the following chain of embeddings.
\begin{align*}
 H^1_2((\delta,\tau); H^{1/2}_2(\Sigma))&\hookrightarrow BU\!C^{1/2}((\delta,\tau); H^{2/2-1/2}_2(\Sigma))\hookrightarrow BU\!C^{1/2}((\delta,\tau); H^{2/q-1/2}_q(\Sigma))\\
&\hookrightarrow BU\!C^{1/2}((\delta,\tau); B^{2/q-1/2}_{qp}(\Sigma)) =BU\!C^{1/2}((\delta,\tau); B^{\mu_c^\kappa-1/p-\kappa/q}_{qp}(\Sigma))
\end{align*}
for any $\tau\in (\delta,t_+)$.
Therefore, $\lim_{t\to t_+} v(t)$ exists in $B^{\mu_c^\kappa-1/p-\kappa/q}_{qp}(\Sigma)$ which according to Theorems 3.2 and  3.3 and parabolic regularization yields global existence for solutions
starting in the critical spaces $X_{crit}= B^{ 2/q -1/2}_{qp}(\Sigma)\times L_q(\Sigma)$, $2< q<12$. These considerations yield the following result on global existence for space dimension $d=3$.

\begin{thm}
Let  $d=3$, $p\geq 2$, $2<q<12$.\\
Then for each initial value ${\sf x}_0=(v_0,w_0)$ and inhomogeneity $\Phi$ with regularity
$$ (v_0,w_0) \in B^{2/q-1/2}_{qp}(\Sigma)\times L_q(\Sigma),\quad  \Phi\in H^1_{p,loc}([0,\infty); W^{-1/q}_q(\Sigma)),$$
the unique solution $(v,w)$ of \eqref{sleveq} exists globally in the class
\begin{equation*}
\begin{split}
&v\in H^2_{p,loc}((0,\infty); W^{-1/q}_q(\Sigma))\cap H^1_{p,loc}((0,\infty); W^{1-1/q}_q(\Sigma))
 \hookrightarrow C^1((0,\infty); B^{1-1/p-1/q}_{qp}(\Sigma)),\\
& w\in H^{3-1/q}_{p,loc}((0,\infty);L_q(\Sigma))\hookrightarrow C^{3-1/p-1/q}((0,\infty);L_q(\Sigma)),
\end{split}
\end{equation*}
and 
$u\in H^1_{p,loc}((0,\infty); H^1_q(\Omega\setminus\Sigma)).$
\end{thm}

The case $q=2$ needs some additional considerations, as Theorem 3.2 does not apply. 
Instead we fix any $\kappa\in[0,1)$ and employ Theorem 3.3.
As above this yields a solution on a maximal interval with regularity
$$ v\in H^2_{p,loc}(0,t_+);H^{-\kappa/2}_2(\Sigma))\cap H^1_{p,loc}((0,t_+); H^{1-\kappa/2}_2(\Sigma)),$$
provided
$$  (v_0,w_0)\in B^{1/2}_{2p}(\Sigma)\times L_2(\Sigma),\quad \Phi\in H^1_{p,loc}([0,\infty);H^{-\kappa/2}_2(\Sigma)).$$
As $ H^{1-\kappa/2}_2(\Sigma))\hookrightarrow L_2(\Sigma)$, the energy estimate \eqref{enbds2} applies, leading to global existence, as before.
 
\begin{corollary} In the situation of Theorem 4.2, let $d=3$, $q=2$, $p\geq2$, $\kappa\in [0,1)$, and asssume
$$ (v_0,w_0)\in B_{2p}^{1/2}(\Sigma)\times L_2(\Sigma),\quad \Phi\in L_{p,loc}([0,\infty);H^{-\kappa/2}_2(\Sigma)).$$
Then the unique solution of \eqref{sleveq} is global and has regularity
\begin{equation*}
\begin{split}
&v\in H^2_{p,loc}((0,\infty); H^{-\kappa/2}_2(\Sigma))\cap H^1_{p,loc}((0,\infty); H^{1-\kappa/2}_2(\Sigma))
 \hookrightarrow C^1((0,\infty); B^{1-1/p-\kappa/2}_{2p}(\Sigma)),\\
& w\in H^{3-\kappa/2}_{p,loc}((0,\infty);L_2(\Sigma))\hookrightarrow C^{3-1/p-\kappa/2}((0,\infty);L_2(\Sigma)),
\end{split}
\end{equation*}
and 
$ u\in L_{p,loc}((0,\infty); H^{1+(1-\kappa)/2}_2(\Omega\setminus\Sigma)).$
These assertions are also valid for $\kappa=1$, $p=2$.
\end{corollary}
\noindent
The last statement in this corollary is proved as follows. For $d=3$ we have the embedding $H^{1/2}_2(\Sigma)\hookrightarrow L_4(\Sigma)$, which by duality implies $L_{4/3}(\Sigma)\hookrightarrow H^{-1/2}_2(\Sigma)$. Next we estimate by H\"older's inequality
$$ |v^2\dot{v}|_{L_{4/3}} \leq |v|_{L_4}^2 |\dot{v}|_{L_4} \leq C |v|_{H^{1/2}_4}^2  |\dot{v}|_{H^{1/2}_2},$$
which implies by \eqref{enbds2} that the nonlinearity $v^2\dot{v}$ in the differentiated equation belongs to $L_{2,loc}((0,\infty);H^{-1/2}_2(\Sigma))$. Therefore, maximal regularity yields $\dot{v}\in H^1_{2,loc}((0,\infty);H^{-1/2}_2(\Sigma))$, provided $\Phi\in L_{2,loc}([0,\infty);H^{-1/2}_2(\Sigma))$.
Then we may employ an approximation argument concerning $\Phi$ to prove the statement.

\section{Stability of Homogeneous Equilibria}
We have already noted that, for $\Phi\equiv0$, problem \eqref{dbe} or \eqref{sleveq} contains the classical FitzHugh-Nagumo system in ODE's as a special case. Therefore, in this case we have homogeneous equilibria, namely ${\sf x}^0 = [0,0]^{\sf T}$ and the nontrivial ones ${\sf x}^\pm =[ v^\pm,w^\pm]^{\sf T}$, where
$$ w^\pm = \frac{c}{b} v^\pm,\quad v^\pm = \frac{a+1}{2} \pm \frac{1}{2} \sqrt{ (a-1)^2 -4\frac{c}{b}}.$$
The nontrivial ones are real if and only if $(a-1)^2\geq 4c/b$. Excluding the degenerate case of a double root, i.e.\ $(a-1)^2= 4c/b$,
it is well-known that for the ODE system, ${\sf x}^0$ is always exponentially stable, ${\sf x}^-$ is a saddle point, hence unstable, and ${\sf x}^+$ is exponentially stable, if and only if
$$ h^\prime(v^+)> (c-b^2)/b,\quad h(v)= v( v^2 -(a+1)v + a+c/b).$$
We want to prove that the same results are valid for \eqref{sleveq} in the framework of $L_q$-theory.

For this purpose, we choose the strong setting as in Section 3.1, with $2d/3q +1/p<5/3$. Then we have the embedding $B^{2-1/p-1/q}_{qp}(\Sigma)\hookrightarrow W^{1+\beta-1/q}_q(\Sigma)$ for $q<d$, which implies $X_{\gamma,1}\hookrightarrow X_\beta$ and allows  to apply the results in Section 5 of \cite{PrSi16} on the qualitative theory of quasilinear parabolic evolution equations, in particular the principle of linear stability; see Theorems 5.3.1, 5.4.1, 5.5.1, 5.6.1, and 5.6.2 in \cite{PrSi16}.

To apply these results, we proceed as in \cite{HiPr18}.
We only have to look at the spectrum of the linearization $L$ of \eqref{sleveq} at the corresponding equilibrium.
For ${\sf x}^0$  this is already clear from  Section 2, to the result that ${\rm Re}\,\sigma(L)\geq \gamma_1>0$, hence this equilibrium is exponentially stable.
For ${\sf x}^\pm$ we find that the linearizations are given by
$$ L^\pm=\left[ \begin{array}{cc} \aA+a+ 3(v^\pm)^2 -2(a+1)v^\pm& 1\\
     -c&b\end{array}\right].$$
We see that $L^\pm=A+B^\pm$, where $B^\pm$ is a relatively compact perturbation of $A$. Such perturbations do not change the essential spectrum, hence we have
$\sigma_{ess}(L^\pm)=\sigma_{ess}(A)=\{b\}$. So we only have to look at the eigenvalues of $L^\pm$.

If ${\sf x}=(v,w)$ is an eigenfunction for eigenvalue $\lambda\neq b$ of $L^\pm$, then
$$ bw=cv+\lambda w,\quad  (\aA+a+ 3(v^\pm)^2 -2(a+1)v^\pm)v +w =\lambda v,$$
which implies
$$ w=\frac{c}{b-\lambda} v,\quad \aA v = \big( \lambda-a -3(v^{\pm})^2+2(a+1)v^\pm +c/(\lambda-b) \big)v,$$
hence $v$ is an eigenvector for $\aA$ with eigenvalue
$$ 0\leq \mu =\lambda-a -3(v^{\pm})^2+2(a+1)v^\pm +c/(\lambda-b).$$
This yields the quadratic equation for $\lambda$
$$ \lambda^2 -{\sf p}^\pm(\mu)\lambda + {\sf q}^\pm(\mu)=0,$$
where
\begin{align*}
{\sf p}^\pm(\mu) &= \mu+a +b + 3(v^{\pm})^2-2(a+1)v^\pm,\\
{\sf q}^\pm(\mu) &= c+b( \mu+a +3(v^{\pm})^2-2(a+1)v^\pm)
\end{align*}
For stability we need ${\sf p}(\mu)$, ${\sf q}(\mu)>0$ for all $\mu\in\sigma(\aA)$, in particular for $\mu=0$, which corresponds to the homogeneous case.
Next we observe that
$$ h^\prime(v^\pm)=3(v^{\pm})^2-2(a+1)v^\pm + a+c/b,$$ hence
$$ {\sf q}^\pm(\mu) = b(\mu+ h^\prime(v^\pm))\quad \mbox{and} \quad {\sf p}^\pm(\mu) = \mu+h^\prime(v^\pm) +b-c/b.$$
It is clear that except for the degenerate case $v^+=v^-$ we have  $h^\prime(v^-)<0$ and $h^\prime(v^+)>0$. Therefore, ${\sf q}^-(0)<0$, 
i.e.\ the equilibrium ${\sf x}^-$ is always unstable, as in the ODE case. On the other hand, if ${\sf p}^+(0)>0$ and ${\sf q}^+(0)>0$, i.e.\ if $h^\prime(v^+)> (c-b^2)/b$, then ${\sf p}^+(\mu),{\sf q}^+(\mu)>0$ for all $\mu\geq0$, hence we have stability of ${\sf x}^+$. This condition is equivalent to stability of ${\sf x}^+$ for the classical FitzHugh-Nagumo system. We also see that equilibrium ${\sf x}^-$ is hyperbolic if and only if
$-h^\prime(v^-) \notin \sigma(\aA)$,  $-h^\prime(v^-)+(c-b^2)/b\notin \sigma(\aA)$ \text{ in case $c-b^2<0$},
and that ${\sf x}^+$ has this property if and only if
$ -h^\prime(v^+)+(c-b^2)/b\not \in \sigma(\aA)$.

This way we have obtained the following result on stability properties of homogeneous equilibria of \eqref{sleveq} or \eqref{dbe}.

\begin{thm} \label{stab-thm}
Let the constants defining $(f,g)$ satisfy~\eqref{abc}, let $d\geq 2$, $p,q\in (1,\infty)$ with $2d/3q +1/p<5/3$, suppose 
$(a-1)^2> 4c/b$, i.e.\ ${\sf x}^-\neq {\sf x}^+$, and let $\Phi\equiv0$.

\smallskip

Then \eqref{sleveq} generates a local semiflow in $X_{\gamma,1}= B^{2-1/p-1/q}_{qp}(\Sigma)\times W^{1-1/q}_q(\Sigma)$. The homogeneous equilibria of this semiflow defined above have the same stability properties as for the classical ODE FitzHugh-Nagumo system, namely 
\begin{itemize}
\item[{\bf (a)}] the equilibrium ${\sf x}^0=(0,0)$ is exponentially stable;
\vspace{1mm}
\item[{\bf (b)}] the equilibrium $ {\sf x}^-$ is unstable. It is hyperbolic if  and only if 
$$ -h^\prime(v^-)\not\in \sigma(\aA),\;  -h^\prime(v^-)+(c-b^2)/b\not\in \sigma(\aA) \text{ in case $c-b^2>0$}; $$ 
\item[{\bf (c)}] the equilibrium ${\sf x}^+$ is stable if and only if $ h^\prime(v^+)>(c-b^2)/b$. It is hyperbolic if and only if
$$-h^\prime(v^+)+(c-b^2)b\not\in \sigma(\aA);$$
\item[{\bf (d)}] If ${\sf x}^\pm$ is hyperbolic, the stable and unstable manifolds exist and are of class $C^1$.
\end{itemize}
\end{thm}
Results like those proved above for the FitzHugh-Nagumo transport can easily be extended to cover a variety of other, more general transport laws. We leave this topic to the interested reader.

\end{document}